\documentstyle{article}%%[12pt]
\def\Bbb{\bf}
\def\kod{\mbox{Kod}}
\def\be{\begin{equation}}
\def\ee{\end{equation}}
\def\bea{\begin{eqnarray*}}
\def\eea{\end{eqnarray*}}
\newtheorem{defn}{Definition}
\newtheorem{lem}{Lemma}
\newtheorem{thm}{Theorem}
\newtheorem{prop}{Proposition}

\newtheorem{conj}{Conjecture}
\newenvironment{xpl}{\mbox{ }\\ \\{\bf Key Example}\mbox{ }}{
\hfill $\Box$\mbox{}\bigskip}

\newenvironment{proof}{\medskip {\bf Proof.}}{\hfill \rule{.5em}{1em} \\}
\begin{document}
\sloppy
\title{On the Scalar Curvature of Complex Surfaces}

\author{Claude LeBrun\thanks{Supported
in part by  NSF grant DMS-9003263.}\\
SUNY Stony
 Brook}

\date{December, 1994}
\maketitle

\begin{abstract} Let $(M,J)$ be a   minimal
 compact complex surface of
K\"ahler type. It is shown that
the smooth 4-manifold $M$ admits a
Riemannian metric of positive scalar curvature
iff $(M,J)$ admits a {\em K\"ahler} metric of
positive scalar curvature. This extends previous results of
Witten and Kronheimer.
 \end{abstract}
%\vfill
%\pagebreak

\bigskip
A {\em complex surface} is
a pair $(M,J)$ consisting of a  smooth compact 4-manifold $M$ and a complex
structure $J$ on $M$;  the latter  means an almost-complex structure tensor
$J: TM\to TM$, $J^2=-1$,  which is locally isomorphic to the usual
constant-coefficient almost-complex
structure on ${\Bbb R}^4={\Bbb C}^2$. Such a complex surface is called
 {\em minimal} if it contains
no embedded copy $C$ of  $S^2$
such that $J(TC ) = TC$ and such that
$C\cdot C= -1$ in homology; this is equivalent to saying that
$(M,J)$ cannot be obtained from another complex surface
$(\check{M},\check{J})$ by the procedure of ``blowing up a point.''

A Riemannian metric $g$ on  $M$ is said to be
is said to be {\em K\"ahler} with respect to $J$ if $g$ is
$J$-invariant  and  $J$ is   parallel with respect to the
metric connection of $g$. If such metrics actually exist,
  $(M,J)$ is then said to be of
 {\em K\"ahler type}; by a
deep result \cite{bpv} of  Kodaira, Todorov, and Siu, this holds
iff  $b_1(M)$ is even.

The purpose of the present note is to prove the following:

\begin{thm}\label{main}
Let $(M,J)$ be a minimal complex surface of K\"ahler type.
Then the following are equivalent:
\begin{description}
\item{(a)}  $M$ admits a Riemannian metric of positive scalar curvature;
\item{(b)} $(M,J)$ admits a K\"ahler metric of positive scalar curvature;
\item{(c)} $(M,J)$ is either   a  ruled surface or ${\Bbb CP}_2$.
\end{description}
\end{thm}
Here a minimal complex surface $(M,J)$ is said to be  {\em ruled} iff  it
is the total space of a holomorphic ${\Bbb CP}_1$-bundle over a compact
complex  curve.

The equivalence between (b) and (c) was proved in one of
Yau's first papers \cite{yau}. By contrast, the link between (a) and (c)
  came to light only  recently,   when
Witten \cite{wit} discovered that a K\"ahler  surface
with $b^+>1$ cannot admit a metric of positive scalar curvature.
Kronheimer \cite{kp} then used a refinement of
 Witten's method to prove that a minimal
surface of general type cannot admit
positive-scalar-curvature metrics. In essence, what will
be shown here
is simply that Kronheimer's method can, with   added care,
 also
be applied to the case of minimal elliptic surfaces.

\section{Seiberg-Witten Invariants}

The   ideas presented  in
 this section are fundamentally due to
Witten \cite{wit}, but much of the formal
framework and many
technical results are due to Kronheimer-Mrowka \cite{KM}.
The work of
Taubes \cite{tau} contains  less elementary but more
robust proofs of other key results presented here. See also \cite{leb}.

Let $(M,g)$ be a smooth compact Riemannian 4-manifold, and suppose that
$M$ admits an almost-complex structure. Then the given component of
the almost-complex structures on $M$ contains almost-complex
structures  $J:TM\to TM$, $J^2=-1$ which are compatible with $g$
in the sense that $J^*g=g$. Fixing such a $J$,  the tangent bundle
$TM$ of $M$ may be given the structure of a rank-2 complex vector bundle
$T^{1,0}$ by defining scalar multiplication by $i$ to be  $J$. Setting
$\wedge^{0,p}:=\wedge^p\overline{T^{1,0}}^*\cong \wedge^pT^{1,0}$,
we may then
 then define rank-2 complex vector bundles
$V_{\pm}$ on $M$ by
\begin{eqnarray} V_+&=& \wedge^{0,0}\oplus \wedge^{0,2}\label{spl}\\
V_-&=&\wedge^{0,1},\nonumber\end{eqnarray}
and $g$ induces canonical Hermitian inner products on these
bundles.

As described, these bundles depend on the choice of a particular
almost-complex structure $J$,
but  they have a deeper meaning \cite{hit} that depends  only on the
homotopy class $c$  of $J$; namely, if we restrict
to a contractible open set  $U\subset M$,  the bundles
$V_{\pm}$ may be canonically
identified with ${\Bbb S}_{\pm}\otimes L^{1/2}$,
where ${\Bbb S}_{\pm}$ are the left- and right-handed
spinor bundles of $g$, and $L^{1/2}$ is a complex
line bundle whose square is the `anti-canonical'
line-bundle $L=(\overline{\wedge^{0,2}})^*\cong \wedge^{0,2}$.
For each   connection $A$ on $L$ compatible with the
$g$-induced inner product, we can thus
define a  corresponding Dirac operator
$$D_{A}: C^{\infty}(V_+)\to C^{\infty}(V_-).$$
If $J$ is parallel with respect to
$g$, so that $(M,g,J)$ is a K\"ahler manifold,
and if $A$ is the Chern connection on the
anti-canonical bundle $L$, then
$D_A ={\sqrt{2}}(\overline{\partial}\oplus \overline{\partial}^*)$,
where  $\overline{\partial}: C^{\infty}(\wedge^{0,0})\to
C^{\infty}(\wedge^{0,1})$ is  the $J$-antilinear part of the exterior
differential $d$, acting on complex-valued functions, and where
$\overline{\partial}^*: C^{\infty}(\wedge^{0,2})\to
C^{\infty}(\wedge^{0,1})$ is the formal  adjoint of
the  map induced by the exterior
differential $d$ acting on 1-forms;
more generally, $D_A$ will differ from
${\sqrt{2}}(\overline{\partial}\oplus \overline{\partial}^*)$
by only $0^{th}$ order terms.

In addition to the metric $g$ and class $c$ of almost-complex structures
$J$, suppose we also choose some  $\varepsilon\in C^{\infty}(\wedge^+)$.
The perturbed Seiberg-Witten equations
\begin{eqnarray} D_{A}\Phi &=&0\label{drc}\\
 iF^+_A+\sigma (\Phi ) &=& \varepsilon \label{qsd}\end{eqnarray}
are then
equations
for an unknown smooth connection $A$ on $L$
and an unknown
smooth section $\Phi$ of $V_+$.
Here the purely imaginary 2-form $F_{A}^+$  is the self-dual part of
the curvature of $A$, and, in terms of (\ref{spl}),
the real-quadratic map $\sigma: V_+\to \wedge^2_+$
is given by
$$\sigma (f, \phi)=(|f|^2-|\phi|^2) \frac{\omega}{4}+ \Im m
(\bar{f}\phi),$$
where $\omega (\cdot, \cdot ) = g(J\cdot, \cdot )$ is the
`K\"ahler' form.

For a fixed metric $g$, let ${\cal M} (g)$ denote the set
of pairs $(A, \Phi)$ which satisfy (\ref{drc}), modulo
the action $(A, \Phi) \mapsto (A + 2d\log u,  u\Phi)$
of the `gauge group'
 of smooth maps   $u: M\to S^1\subset {\Bbb C}$.
We may then view (\ref{qsd}) as defining
 a map $\wp : {\cal M} (g)
\to C^{\infty}(\wedge^+)$. One can show \cite{KM} that this
is a proper map, and so, in particular, has compact fibers.

A solution
$(A , \Phi)$
is called {\em reducible} if $\Phi\equiv 0$; otherwise, it is
called {\em irreducible}. Let ${\cal M}^* (g)$ denote the
image in ${\cal M} (g)$ of the
set of irreducible solutions. Then  ${\cal M} (g)$
is \cite{KM} a smooth Fr\'echet manifold, and an
index calculation shows that the smooth map
$\wp: {\cal M}^* (g)
\to C^{\infty}(\wedge^+)$ is generically finite-to-one.
 Let us define a solution  $(\Phi , A)$  to be {\em transverse} if
it corresponds to a regular point of $\wp$. This holds iff
 the linearization
$C^{\infty}(V_+\oplus \wedge^1)\to C^{\infty}(\wedge_+^2)$
of the left-hand-side of (\ref{qsd}), constrained by the
linearization of
(\ref{drc}), is  surjective.

\begin{xpl}
Let $(M,g,J)$ be a K\"ahler surface, and let    $s$
denote the scalar curvature of $g$. Set
   $\varepsilon = (s+1)\omega /4$,
 set $\Phi = (1, 0)\in
\wedge^{0,0}\oplus \wedge^{0,2}$, and let $A$ be
the  Chern connection on the anti-canonical bundle.
Since $iF_A$ is just the Ricci form of $(M,g,J)$, it follows that
 $iF^+_A+\sigma ( \Phi )=s\omega /4 +\omega/4=\varepsilon$,
and $(\Phi , A)$ is thus
an irreducible solution of  equations
(\ref{drc}) and (\ref{qsd}).

The linearization of (\ref{drc}) at this solution
is just
\be (\overline{\partial}\oplus \overline{\partial}^*)(u+\psi)=
-\frac{1}{2}\alpha ,\label{subs}\ee
where $(v,\psi)\in C^{\infty}(V_+)$
is the linearization of $\Phi=(f,\phi)$
and
$\alpha \in \wedge^{0,1}$ is
the $(0,1)$-part of the purely imaginary 1-form which is the
linearization of $A$. Linearizing
 (\ref{qsd}) at our solution yields the operator
$$(v,\psi, \alpha) \mapsto
id^+(\alpha-\bar{\alpha})+\frac{1}{2}(\Re e v )\omega
+\Im m \psi . $$
Since the right-hand-side is a real self-dual form, it
is completely characterized by its component in the
$\omega$ direction and its $(0,2)$-part.
The $\omega$-component of this operator is just
$$  (v,\psi, \alpha)\mapsto \Re e \left[-\bar{\partial}^*\alpha
+ \frac{v}{2} \right], $$ while the $(0,2)$-component
is $$(v,\psi, \alpha)\mapsto i\bar{\partial} \alpha - \frac{i}{2}
\psi.$$
Substituting (\ref{subs}) into these expressions,
we obtain the operator
\bea C^{\infty}({{\Bbb C}}\oplus \wedge^{0,2})&\longrightarrow
 &C^{\infty}(
{\Bbb R}\oplus \wedge^{0,2})\\
(v,\psi)~~~~&\mapsto & ( \Re e
\left[\Delta +\frac{1}{2}\right] v, -i
\left[\Delta +\frac{1}{2}\right]\psi)  ,
\eea
which is surjective because  $\Delta  + \frac{1}{2}$ is a
positive self-adjoint elliptic operator.
 The constructed
solution is therefore transverse.
\end{xpl}

For any metric $g$, let $c_1^+$ denote the image
of $c_1(L)\in H^2(M,{\Bbb R})$ under
orthogonal\footnote{with respect to the intersection form} projection
 to the linear subspace $H^+(g)$
of deRham classes which are represented by self-dual 2-forms with
respect to $g$. Given any $\varepsilon\in C^{\infty}(\wedge^+)$, let
$\varepsilon_H$ denote its harmonic part; this is a closed
self-dual 2-form, since the   Laplacian
commutes with the Hodge star operator.

\begin{lem} Let $g$ be any Riemannian metric on $M$, and let
$\varepsilon\in C^{\infty}(\wedge^+)$.
Suppose that $[\varepsilon_H]\neq 2\pi c_1^+$ in deRham cohomology.
Then every solution of (\ref{drc}) and (\ref{qsd}) is
irreducible.
\end{lem}
\begin{proof}
If $\Phi \equiv 0$, (\ref{qsd}) says that
$c_1(L)$ is represented by $\varepsilon/2\pi$ plus
an anti-self-dual form. Taking the harmonic part of this
representative and projecting to the self-dual harmonic forms
then yields $c^+_1=[\varepsilon_H]/2\pi$.
\end{proof}

\begin{defn}
If $g$ is a smooth Riemannian metric on $M$ and
$\varepsilon\in C^{\infty}(\wedge^+)$ is such that
$[\varepsilon_H]\neq 2\pi c_1^+$, then, with respect to
$c=[J]$, we will say that $(g,\varepsilon)$ is a {\em good pair}.
The path components of the manifold of all good pairs $(g,\varepsilon)$
will
be called {\em chambers}.
\end{defn}

\begin{lem}\label{wind}
If $b^+(M)> 1$, there is exactly one chamber. If
$b^+(M)=1$, there are exactly two chambers.
\end{lem}
\begin{proof}
The projection $(g, \varepsilon)\mapsto g$ factors through
the rank-$b^+$ vector  bundle $H^+$ over the space of
Riemannian metrics via a map
$(g, \varepsilon)\mapsto (g, [\varepsilon_H])$ with
connected fibers.
Now
$2\pi c^+_1(g)$ is a smooth section of
 $H^+$, and   since the space of Riemannian metrics
is path-connected,  the
  complement of this section is connected if $b^+>1$, and
has exactly two components if $b^+=1$. The result follows.
\end{proof}

\begin{lem}\label{kind}
Suppose that $b^+(M)=1$, $c_1(L)\neq 0$, and $c_1^2(L)\geq 0$.
Then $(g,0)$ is a good pair  for any metric $g$. In particular,
the chamber containing $(g,0)$ is independent of $g$, and will be
called the {\em preferred chamber}.
\end{lem}
\begin{proof}
Our hypotheses say the intersection form on $H^2$ is a Lorentzian
inner product and that
 $c_1(L)$ is a non-zero null or time-like vector. Since
 $2\pi c_1^+(g)$ is the projection of $2 \pi c_1(L)$ onto
a time-like line, it thus never equals
$0=[0_H]$.
\end{proof}

\begin{defn}
Let $(M,c)$ be a compact 4-manifold equipped with
a  homotopy class $c=[J]$ of almost-complex structures.
A good pair $(g, \varepsilon)$ will be called {\em excellent}
if $\varepsilon$ is a regular value of the map
$\wp : {\cal M}^*(g)\to C^{\infty}(\wedge^+)$.
\end{defn}

\noindent Notice that $(g, \varepsilon )$ is excellent
iff every solution of (\ref{drc}) and (\ref{qsd}) with
respect to  $(g, \varepsilon )$ is irreducible and transverse.

\begin{defn}\label{inv}
Let $(M,c)$ be a compact 4-manifold equipped with
a   class $c=[J]$ of almost-complex structures.
 If
$(g,\varepsilon)$ is an excellent pair on $M$, we define
its (mod 2) {\em Seiberg-Witten invariant} $n_c(M,g,\varepsilon)
\in {\Bbb Z}_2$
to be
$$n_c(M,g, \varepsilon)=\# \{\mbox{gauge classes of
solutions of  (\ref{drc})   and  (\ref{qsd})}
\}
\bmod 2 $$
calculated with respect to $(g, \varepsilon)$.
\end{defn}

\begin{lem} If two excellent pairs are in the same chamber,
they have the same Seiberg-Witten invariant $n_c$.
\end{lem}
\begin{proof} Any two such pairs can be joined by a
path of good pairs which is transverse to $\wp$. This gives a
cobordism between the relevant solution spaces.
\end{proof}

\begin{defn}
Let $(M, c)$ be a smooth 4-manifold equipped with
a class of almost-complex structures. If $b^+(M) > 1$,
the Seiberg-Witten invariant $n_c(M)$ is defined to be
$n_c(M,g,\varepsilon)$, where $(g,\varepsilon)$ is
any excellent pair. If $b^+(M) = 1$,
$c_1 (L)\neq 0$, and $c_1^2(L)\geq 0$, then
the Seiberg-Witten invariant $n_c(M)$ is defined to be
$n_c(M,g,\varepsilon)$, where $(g,\varepsilon)$ is
any excellent pair in the preferred chamber.
\end{defn}

\begin{thm}\label{non}
Let $(M,J)$ be a compact complex surface which
admits a
K\"ahler metric $g$; let $c=[J]$.
 Then there is a chamber for which
the Seiberg-Witten invariant $n_c$ is non-zero.
Moreover, if   $c_1\cdot [\omega]<0$,
where  $[\omega]$ is the K\"ahler class of $g$,
then the chamber in question  contains the good pair $(g,0)$.
 \end{thm}
\begin{proof} Set $\varepsilon= (s+1)\omega/4$, where
$s$ is the scalar curvature of the K\"ahler metric
  $g$. Then $\varepsilon_H= (s_0+1)\omega/4$, where
the average value $s_0$ of the scalar curvature
of $g$ is given by
$$
s_0=\frac{\int_M s ~d\mu}{\int_M d\mu}=
\frac{2\int_M \rho\wedge\omega}{\int_M \omega\wedge\omega/2} =
 8\pi
\frac{ c_1\cdot [\omega]}{[\omega]^2}
$$
because the Ricci form $\rho$ represents $2\pi c_1$.
On the other hand, $2\pi c_1^+$ is represented by
the harmonic form $s_0\omega/4$, so we always have
$2\pi c^+_1\neq [\varepsilon_H]$. This shows that $(g, \varepsilon )$
is a good pair. Moreover, if $c_1\cdot [\omega ] < 0$, then
$s_0 < 0$,  and
and $(g, t\varepsilon )$ is a good pair for all
$t\in [0,1]$. Thus it suffices to show that
$n_c(M,g,\varepsilon )\equiv 1 \bmod 2$.

We will accomplish this by showing  that, with
respect to
$(g, \varepsilon )$ and
up to gauge equivalence, there is exactly one solution
of the perturbed Seiberg-Witten
equations, namely the transverse solution  described in the
Key Example. Indeed, suppose that $\Phi = (f,\phi )$
is any solution of (\ref{drc}) and (\ref{qsd}), and let
$\hat{\Phi}=(f,-\phi)$.
The  Weitzenb\"ock formula for the twisted Dirac operator
and equation (\ref{qsd}) thus tell us that
\bea 0=D_A^*D_A\Phi&=&\nabla^*\nabla\Phi+ \frac{s}{4}\Phi
+ \frac{1}{2} F_A\cdot \Phi\\
&=&\nabla^*\nabla\Phi+ \frac{s}{4}\Phi
+ \frac{i}{2}\sigma (\Phi)\cdot \Phi
- \frac{i}{2}\varepsilon\cdot \hat{\Phi}
\\&=&\nabla^*\nabla\Phi+ \frac{1}{4}(s+|\Phi|^2)\Phi
-  \frac{1}{4}(s+1) \hat{\Phi}
\eea
because the $\pm 1$-eigenspaces of Clifford multiplication
on $V_+$
by $-2i\omega$ are respectively $\wedge^{0,0}$ and $\wedge^{0,2}$.
  Projecting to the  $\nabla$-invariant sub-bundle
$\wedge^{0,0}\subset V_+$ now yields
$$0 =4\nabla^*\nabla f - f + |\Phi|^2f .$$
At the maximum of $|f|^2$, we thus have
\bea
0\leq 2\Delta \langle f, f \rangle &= & 4\langle \Delta f, f \rangle -
4\langle \nabla f, \nabla f \rangle
\\&\leq & |f|^2 - |\Phi|^2|f|^2
\\&\leq & (1 - |f|^2)  |f|^2\eea
because $|\Phi|^2 = |f|^2 + |\phi |^2\geq |f|^2$. This
gives us the inequality
\be
|f|^2\leq 1 \label{est}\ee
at all points of $M$, with equality only
when $\phi=0$ and $\nabla f = 0$.

On the other hand,the closed 2-form $iF_A$ is  in the
same cohomology class as the Ricci form $\rho$, which
satisfies $\langle \omega , \rho \rangle = s/2$.
Writing $iF_A=\rho + d \beta$ for some 1-form $\beta$, we have
\bea  \int_M \left(|f|^2-|\phi|^2\right) d\mu
&= & 2\int_M \langle \omega , \sigma (f,\phi ) \rangle d\mu
\\&=& 2\int_M \langle \omega , -iF_A+ \varepsilon \rangle d\mu
\\&=& 2\int_M \langle \omega , -\rho-d\beta\rangle d\mu +
\int_M  \langle \omega ,   (s+1)\frac{\omega}{2} \rangle d\mu
\\&=&-\int_M    s   ~ d\mu -
2\int_M   \langle d^*\omega , \beta\rangle  d\mu +
\int_M    (s+1) d\mu \\&=&
 \int_M  1~ d\mu ,\eea
which is to say  that
$$\int_M (|f|^2 -1)d\mu =\int_M |\phi|^2 d\mu \geq 0.$$
The $C^0$ estimate (\ref{est}) thus implies that $|f|^2\equiv 1$,
 $\phi\equiv 0$, and
$\nabla f\equiv 0$. The connection $\nabla$ induced on
the $\wedge^{0,0}$ by $A$ is therefore
flat and trivial, and $A$ is thus
gauge equivalent to the Chern connection on $L$. Hence
 our solution    coincides, up to gauge
transformation, with that of the example; in particular, every solution
with respect to $(g,\varepsilon)$ is irreducible and transverse,
and
$(g,\varepsilon)$ is an excellent pair.  But since there is only
one gauge class of
solutions with respect to $(g,\varepsilon)$, it follows that
 $n_c =1\bmod 2$ for the chamber containing $(g, \varepsilon )$.
 \end{proof}

\begin{thm}\label{pos}
Let $M$ be a compact 4-manifold which admits a class
$c=[J]$ of almost-complex structures and has $b^+  > 0$.
If $g$ is a metric of positive scalar curvature on $M$, then
$(g,0)$ is in the closure of a  chamber for which
   $n_c =0$.
\end{thm}
\begin{proof} Suppose not. For every $\epsilon > 0$,
there is a    $\varepsilon$ such that
 $\sup |\varepsilon | < 2\epsilon$ and such that
$(g, \varepsilon)$ is an excellent pair.
If $n_c(M,g,\varepsilon )\neq 0$, there is a
solution $\Phi\not\equiv 0$ of
equations (\ref{drc}) and (\ref{qsd}) with respect to
$g$ and  $\varepsilon = 0$. The Weitzenb\"ock formula
$$0=D_A^*D_A\Phi = \nabla^*\nabla \Phi + \frac{s+|\Phi|^2}{4}\Phi
-\frac{i}{2}\varepsilon \cdot \Phi
$$
then implies  that
$$0 > \int_M\left( \frac{s - \epsilon}{4} \right)|\Phi|^2 d\mu .
$$
Taking $\epsilon < \min s$ then yields a contradiction.
\end{proof}

\section{Surface Classification and Scalar Curvature}

Recall \cite{bpv}
that the Kodaira dimension $\kod (M,J) \in \{ -\infty , 0, 1, 2\}$
of a compact
complex surface $(M,J)$ is defined to be $\lim\sup (\log h^0(M, {\cal O}
(L^{*\otimes m}))/\log m )$. The following well-known degree argument
may be found e.g. in
\cite{yau}.

\begin{lem}  Let $[\omega ]$ be a K\"ahler class on
a compact complex surface $(M,J)$ of $\kod \geq 0$. Then
$c_1\cdot [\omega ]\leq 0$, with
equality iff $(M,J)$ is a minimal surface of
$\kod = 0$.
\end{lem}
\begin{proof}
If $\kod (M,J) \geq 0$, some positive
  power $\kappa^m$ of the canonical
line bundle has a holomorphic section. Let
${\bf D}$ be the holomorphic curve, counted with appropriate
multiplicity, where this section vanishes.
The homology class $[{\bf D}]\in H_2(M)$
is then the Poincar\'e dual of $c_1(\kappa^m)=-mc_1(L)$.
The area of ${\bf D}$ is thus
$$\int_{\bf D}\omega = -mc_1\cdot [\omega]$$
which shows that $c_1\cdot [\omega ] \leq 0$,
with equality iff ${\bf D}=\emptyset$. Since the latter
happens iff  $\kappa^m$ is holomorphically trivial,
the result follows.
\end{proof}

 This leads us directly to   a
result  first discovered by Kronheimer
\cite{kp}.

\begin{prop}(Kronheimer)
  Let $(M,J)$ be a  minimal complex surface of
$\kod = 2$. Then   $M$ does not admit a Riemannian metric
of positive scalar curvature.
\end{prop}
\begin{proof}
Such a surface is automatically \cite{bpv} of K\"ahler type,
and has $c_1^2 > 0$. The Seiberg-Witten invariant $n_c(M)$ is
of $(M,J)$ is thus well-defined by Lemma \ref{kind}
and is non-zero by virtue of Theorem \ref{non}. The
result therefore follows by Theorem \ref{pos}.
\end{proof}

 Similar reasoning yields

\begin{prop}
  Let $(M,J)$ be a  minimal complex surface of K\"ahler type with
$\kod = 1$. Then   $M$ does not admit a Riemannian metric
of positive scalar curvature.
\end{prop}
\begin{proof}
Such a surface  must \cite{bpv}  have
$c_1^2 = 0$ and $c_1\neq 0$. The Seiberg-Witten invariant $n_c(M)$
of $(M,J)$ is thus well-defined by Lemma \ref{kind}, and the
conclusion now follows by the same argument used above.
\end{proof}

The next case is actually covered by existing
results \cite{lic,gl,sy}, but  a Seiberg-Witten proof is
given for the sake of completeness.

\begin{prop}  Let $(M,J)$ be a minimal complex surface of K\"ahler
type such that
$\kod (M,J) =0$.  Then $M$ does not admit a Riemannian metric
of positive scalar curvature.
\end{prop}
\begin{proof} Any such an $M$ is finitely
covered by a surface $\tilde{M}$
with $b^+=3$; in fact,  $\tilde{M}$
is either a K3 surface or a
  4-torus.
The Seiberg-Witten invariant $n_c(\tilde{M})$
is thus well-defined by Lemma \ref{wind}
and is non-zero by virtue of Theorem \ref{non}.
By Theorem \ref{pos}, $\tilde{M}$ does not
admit a metric with $s>0$. The result therefore follows because
any metric on $M$ can be pulled back to $\tilde{M}$.
\end{proof}

 Our next result immediately implies  Theorem \ref{main}:

\begin{thm}
Let $(M,J)$ be a minimal surface of
K\"ahler type. If $M$ admits a Riemannian metric of
positive scalar curvature, then $(M,J)$
is either ${\Bbb CP}_2$ or a ruled surface.
As a consequence, $(M,J)$ therefore carries {\em K\"ahler}
metrics of positive scalar curvature.
\end{thm}
\begin{proof}
By the proceeding Propositions, $(M,J)$ must have Kodaira
dimension $-\infty$. The Kodaira-Enriques classification  \cite{bpv}
thus says that $(M,J)$
is either ${\Bbb CP}_2$ or a ruled surface. Now \cite{yau}
any minimal ruled surface admits K\"ahler metrics of
positive scalar curvature; indeed, if $M\cong {\Bbb P}(E)$,
where
 $\varpi : E\to C$ is rank-2 holomorphic vector bundle, then,
for   any K\"ahler form $\omega_C$ on the
Riemann surface $C$ and any   Hermitian norm $h:E\to {\Bbb R}$
  on the complex vector bundle $E$, the $(1,1)$-form
$$\omega = \varpi^*\omega_C + \epsilon\left(i\partial\bar{\partial}
\log h\right)$$
is a K\"ahler form on $M$ with positive scalar curvature
if $\epsilon > 0$ is sufficiently small.
Since the Fubini-Study metric
on ${\Bbb CP}_2$ is also a K\"ahler metric
 of
positive scalar curvature, the result follows.
\end{proof}

Since the minimality hypothesis only features as
a technicality in connection with the $b^+=1$ case,
the following conjecture now seems extremely credible:

\begin{conj}
Let $(M,J)$ be a compact  complex surface with $b_1$ even.
Then the following are equivalent:
\begin{description}
\item{(a)}  $M$ admits a Riemannian metric of positive scalar curvature;
\item{(b)} $(M,J)$ admits a K\"ahler metric of positive scalar curvature;
\item{(c)} $(M,J)$ is either ${\Bbb CP}_2$ or a blow-up of
some  minimal ruled surface.
\end{description}
\end{conj}

However,  even (b) $\Leftrightarrow$ (c) is
only known `generically;' cf. \cite{hit2,klp,ls}.


\begin{thebibliography}{99}

 \bibitem{bpv}  W. Barth, C. Peters, and A. Van de Ven,
{\bf Compact Complex Surfaces}, Springer-Verlag, 1984.

 \bibitem{gl} M. Gromov and H.B. Lawson,
``Spin and Scalar Curvature in teh Presence of the
Fundamental Group,'' {\bf Ann.\ Math.\ 111} (1980) 209--230.


\bibitem{hit} N.J. Hitchin, {\em Harmonic Spinors}, {\bf Adv.\ Math.\ 14}
(1974) 1--55.


\bibitem{hit2} N.J. Hitchin, {\em On the Curvature of
Rational Surfaces}, {\bf Proc. Symp. Pure Math. 27} (1975) 65--80.

\bibitem{klp}
 J.-S. Kim, C. LeBrun  and M. Pontecorvo, {\em
Scalar-Flat K\"ahler Surfaces of All Genera},
 preprint, 1994.


\bibitem{kp}  P. Kronheimer, private communication.

 \bibitem{KM}  P. Kronheimer and T. Mrowka,
{\em The Genus of Embedded Surfacs in the
Complex Projective Plane}, {\bf Math.\ Res.\ Lett.\ }
to appear.


 \bibitem{leb}  C. LeBrun, {\em
Einstein Metrics and Mostow Rigidity},
{\bf Math.\ Res.\ Lett.\ }
to appear.


 \bibitem{ls}  C. LeBrun and S. Simanca,
 {\em  On K\"ahler Surfaces of Constant
Positive Scalar Curvature},
{\bf J. Geom.\ Analysis,} to appear.




\bibitem{lic} A. Lichnerowicz,
{\em Spineurs Harmoniques}, {\bf C.R. Acad. Sci. Paris 257}
(1963) 7--9.


\bibitem{sy}  R. Schoen and S.-T. Yau, {\em The Structure
of Manifolds of Positive Scalar Curvature,}
 {\bf Man.\  Math.\ 28} (1979) 159--183.




 \bibitem{tau}
C.H. Taubes, {\em The Seiberg-Witten Invariants and
Symplectic Forms}, {\bf Math.\ Res.\ Lett.\ }
to appear.

\bibitem{wit} E. Witten, {\em Monopoles and Four-Manifolds},
preprint, 1994.


\bibitem{yau}  S.-T. Yau, {\em On the Curvature of Compact
Hermitian Manifolds,} {\bf Inv.\  Math.\ 25} (1974) 213-239.

\end{thebibliography}
\end{document}